\documentclass[]{interact}

\usepackage{epstopdf}
\usepackage{subfigure}

\usepackage[numbers,sort&compress]{natbib}
\bibpunct[, ]{[}{]}{,}{n}{,}{,}

\theoremstyle{plain}

\theoremstyle{definition}

\theoremstyle{remark}

\begin{document}

\articletype{ARTICLE TEMPLATE}

\title{An efficient method to compute solitary wave solutions of fractional Korteweg-de Vries equations}

\author{
\name{A. Dur\'an\thanks{CONTACT A. Dur\'an. Email: angel@mac.uva.es}}
\affil{Department of Applied Mathematics, University of Valladolid, Paseo Bel\'en, 15, 47011, Valladolid, Spain.}
}

\maketitle

\begin{abstract}
Considered here is an efficient technique to compute approximate profiles of solitary wave solutions of fractional Korteweg-de Vries equations. The numerical method is based on a fixed-point iterative algorithm along with extrapolation techniques of acceleration. This combination improves the performance in both the velocity of convergence and the computation of profiles for limiting values of the fractional parameter. The algorithm is described and numerical experiments of validation are presented. The accuracy attained by the procedure can be used to investigate additional properties of the waves. This approach is illustrated here by analyzing the speed-amplitude relation.
\end{abstract}


\section{Introduction}

The paper is concerned with the computation of solitary wave solutions of the fractional Korteweg-de Vries (fKdV) equation
\begin{equation}
u_{t}+u^{p}u_{x}-(D^{\alpha}u)_{x}=0.\label{ijcm1}
\end{equation}
In (\ref{ijcm1}), $u=u(x,t)$ is a real-valued function of $x\in\mathbb{R}, t\geq 0$, $p\in\mathbb{N}, \alpha\in\mathbb{R}$ and $D^{\alpha}$ stands for the linear operator represented
by the symbol
\begin{equation}
\widehat{\left(D^{\alpha}g\right)}(\xi)=\beta(\xi)\widehat{g}(\xi),\quad \beta(\xi)=|\xi|^{\alpha}, \label{ijcm2}
\end{equation}
where
\begin{equation*}
\widehat{g}(\xi)=\int_{-\infty}^{\infty}g(x)e^{-i\xi x}dx,\quad \xi\in\mathbb{R}
\end{equation*}
is the Fourier transform, defined on the space of squared integrable functions $g\in L^{2}(\mathbb{R})$. Equation (\ref{ijcm1}) is relevant as a dispersive and nonlinbear perturbation of Burgers' inviscid equation. The parameters $\alpha$ and $p$ govern, respectively, the dispersion and the nonlinear effects and in this sense (\ref{ijcm1}) is suitable to investigate the relations between nonlinearity and dispersion that lead to different dynamics, such as existence and stability of solitary waves, blow-up phenomena, etc, \cite{KleinS2015}.
The case $\alpha=2$ corresponds to the classical Korteweg-de Vries (KdV, when $p=1$) and generalized Korteweg-de Vries (gKdV, when $p\geq  2$) equations, while $\alpha=1$ leads to the Benjamin-Ono (BO) equation and generalized (gBO) versions.

The parameters $\alpha$ and $p$ also determine several mathematical properties of the initial value problem for (\ref{ijcm1}), (\ref{ijcm2}). The main results on well-posedness in the literature concern the case $p=1$, for which the Cauchy problem, when $\alpha\geq 1$, is proved to have global solutions in suitable functional spaces, \cite{HerrIKK2010}. When $-1<\alpha<0$, suitable smooth initial data lead the corresponding solution to blow up at finite time, \cite{KleinS2015}. As far as the case $0<\alpha<1$ is concerned, see \cite{LinaresPS2014} for local well-posedness results. Global weak solutions, without uniqueness, in
\begin{equation*}
L^{\infty}(\mathbb{R},H^{\alpha/2}(\mathbb{R}))=\{v:\mathbb{R}\rightarrow H^{\alpha/2}(\mathbb{R})/\max_{t\in\mathbb{R}}||v(t)||_{H^{\alpha/2}(\mathbb{R})}<\infty\},
\end{equation*}
for initial data in the Sobolev space $H^{\alpha/2}(\mathbb{R})$ (with the usual norm $||\cdot||_{H^{\alpha/2}(\mathbb{R})}$, are proved to exists in \cite{Saut1979} (see also \cite{LinaresPS2014}) for $\alpha>1/2$. Also, Klein and Saut, \cite{KleinS2015} conjecture several cases of blow-up with different structure: first, no hyperbolic blow-up (blow-up of the spatial gradient with bounded sup-norm) exists; when $1/2<\alpha<1$, the solution is global; when $1/3<\alpha<1/2$, there is a sort of nonlinear dispersive blow-up, \cite{MartelMR2014I,MartelMR2014II,MartelMR2014III}, and, finally, blow-up of different type occurs when $0<\alpha<1/3$.

On the other hand, at least formally, the following quantities are preserved by smooth enough, decaying solutions of (\ref{ijcm1})
\begin{equation*}
C(u)=\int_{-\infty}^{\infty} u(x,t)dx,\label{ijcm3a}
\end{equation*}
\begin{equation*}
M(u)=\int_{-\infty}^{\infty} u^{2}(x,t)dx,\label{ijcm3b}
\end{equation*}
\begin{equation}
M(u)=\int_{-\infty}^{\infty} \left(\frac{1}{2}|D^{\alpha/2}u(x,t)|^{2}-
\frac{1}{(p+1)(p+2)}u^{p+2}(x,t)\right)dx.\label{ijcm3c}
\end{equation}
The quantity (\ref{ijcm3c}) is well defined when $\alpha\geq 1/3$ and provides a Hamiltonian structure to (\ref{ijcm1}), see \cite{KleinS2015} for the case $p=1$ and \cite{Angulo2017} for $p>1$.

An additional, relevant point on the dynamics of (\ref{ijcm1}) concerns the existence and stability of solitary wave solutions. They are solutions of the form $u(x,t)=\phi(x-ct)$ for some speed $c>0$ and profile $\phi_{c}=\phi_{c}(X)$ with $\phi_{c}(X)\rightarrow 0$ as $|X|\rightarrow\infty$. Substituting into (\ref{ijcm1}) and integrating once, the profile $\phi_{c}$ must satisfy
\begin{equation}
D^{\alpha}\phi_{c}+c\phi_{c}-\frac{\phi_{c}^{p+1}}{p+1}=0.
\label{ijcm4}
\end{equation}
Some known results on existence and stability of solutions of (\ref{ijcm4}) are now summarized.
For the case $p=1$, solitary waves are proved to exist when $\alpha>1/3$, \cite{EhrnstromGW2012,FelmerQT2012,FrankL2013}, with an asymptotic decay as $1/x^{1+\alpha}, |x|\rightarrow\infty$; orbital stability, \cite{Angulo2009,GrillakisSS1987},  for $\alpha>1/2$ is proved in \cite{LinaresPS2014}, while some spectral instability analysis can be obtained from \cite{KapitulaS2014,Pelinovsky2013,Angulo2017}. In the case $p>1$, as mentioned in \cite{Angulo2017} (see also references therein), existence of solutions of (\ref{ijcm4}) holds for $\alpha>1$ and any $p$, with orbital stability when $p<2\alpha$. For $0<\alpha<1$, existence can be derived for $1<p<2\alpha/(1-\alpha)$, \cite{Weinstein1986,Weinstein1987,Arnesen2016}. Results on orbital stability when $1/2<\alpha<1, p<2\alpha$, as well as a linear instability criterium, can be seen in \cite{Angulo2017}.

Since no explicit formulas for solitary wave solutions of (\ref{ijcm1}) are, except in the classical cases $\alpha=1,2$, unknown, then some numerical method for the generation of approximate profiles is required. In this sense, Klein and Saut, \cite{KleinS2015}, solve numerically (\ref{ijcm4}) for $p=1$ to construct approximate solitary wave profiles. The numerical method to this end is based on transforming (\ref{ijcm4}) into the corresponding algebraic equations for the Fourier transform of the profile, which are iteratively solved by Newton iteration. The procedure is implemented by approximating each profile with a trigonometric interpolant polynomial on a long enough interval, in such a way that the Newton method is applied to the system for the corresponding Fourier components. As mentioned in \cite{KleinS2015}, the algebraic decrease of the modulus of the Fourier coefficients, due to the loss of smoothness of the periodic approximations of the profiles at the boundary of the computational domain (the solitary wave profile decreases slowly at infinity) leads to a slow convergence of the iteration. This is addressed by taking a large computational domain and a high number of Fourier modes. The resolution is also performed by using GMRES, \cite{SaadS1986}, to compute iteratively the inverse of the Jacobian matrix.

In this paper, an alternative to construct numerically solitary wave profiles of (\ref{ijcm4}) is proposed. The technique was successfully applied to the numerical generation of periodic traveling wave solutions of the fKdV equation in \cite{AlvarezD2017}. The main points of this approach here are the following:
\begin{itemize}
\item The method is also based on the implementation of (\ref{ijcm4}) in Fourier space for the periodic approximation on a long enough interval.
\item The algebraic system for the discrete Fourier coefficients of the trigonometric interpolant is however iteratively solved by the Petviashvili method, \cite{Pet1976}, a fixed point type algorithm which may overcome some of the limitations of the Newton iteration.
\item In order to improve the slow convergence due to the periodic approximation, the Petviashvili method is complemented with the use of acceleration techniques based on extrapolation, \cite{Sidi2003,sidifs}, which have shown a relevant performance in the numerical generation of solitary waves, \cite{AlvarezD2016}.
 \end{itemize}  
 The main contributions are the following:
 \begin{itemize}
 \item The combination of the Petviashvili method with extrapolation improves the computation of the solitary wave profiles. The improvement is observed in mainly two points: the first one is that the method is able to generate numerical profiles for limiting values of $\alpha$. A second point of improvement is found in the efficiency, since with a high number of Fourier modes and on a long interval, the extrapolation technique accelerates the convergence in a relevant way.
 \item These new advantages can be used to study computationally additional properties of the waves, such as the speed-amplitude relation.
 \item The method can also be applied to compute approximate solitary wave solutions of other generalizations of (\ref{ijcm1}) of the form, \cite{Angulo2017}
 \begin{equation}
 u_{t}+(f(u))_{x}-({\cal M}u)_{x}=0,\label{ijcm5}
 \end{equation}
 where ${\cal M}$ is a linear, pseudo-differential operator associated to a continuous, even, real-valued Fourier symbol $\beta(\xi)$ and $f$ is a smooth, nonlinear, real-valued function. A relevant example of (\ref{ijcm5}) is the extended Whitham equation, \cite{Lannes2013,LannesS2013,DinvayMDK2017}, for which
 \begin{equation}
 \beta(\xi)=\left(1+\gamma|\xi|^{2}\right)^{1/2}\left(\frac{\tanh{\xi}}{\xi}\right)^{1/2},\quad f(u)=\frac{u^{2}}{2}.\label{ijcm6}
 \end{equation}
 In (\ref{ijcm6}), the parameter $\gamma \geq 0$ controls the surface tension effects in the model. The case $\gamma=0$ leads to the classical Whitham equation, \cite{Whitham1967,MoldabayevKD2015}. The computation of traveling-wave solutions of equations of the form (\ref{ijcm5}), kncluding the extended Whitham equation (\ref{ijcm6}), has been made in the literature with different techniques, see e.~g. the references in \cite{AlvarezD2017} and, more recently, the method introduced and performed in \cite{KalischMV2017}, based on continuation with spectral projection. 

The technique can also be applied to study the solitary wave solutions of fractional BBM type equations
 \begin{equation*}
 u_{t}+u_{x}+(f(u))_{x}+({\cal M}u)_{t}=0,
 \end{equation*}
 with ${\cal M}$ and $f$ defined as in (\ref{ijcm5}), see \cite{Angulo2017}.
  \end{itemize} 
  The structure of the paper is as follows. In Section \ref{sec2} the numerical method to compute approximate solitary profiles, based on the Petviashvili method and accelerating techniques, is described, along with some implementations details. The purpose of Section \ref{sec3} is two-fold: a first group of experiments validates the efficiency of the method and studies its performance. This is used, in a second part, to analyze computationally additional properties of the solitary waves. The illustration is focused on the speed-amplitude relation and its dependence on the parameters $\alpha$ and $p$.

\section{Numerical generation of solitary waves}
\label{sec2}
\subsection{The Petviashvili method}
In order to describe the numerical method to compute solitary waves of (\ref{ijcm1}), observe that (\ref{ijcm4}) can be written in the form
\begin{equation}\label{ijcm21}
\mathcal{L}\phi=\mathcal{N}(\phi),\quad \mathcal{L}=D^{\alpha}+c,\quad
\mathcal{N}(\phi)=\frac{\phi^{p+1}}{p+1},\quad c>0.
\end{equation}
Note now that the nonlinear term $\mathcal{N}$ is homogeneous of degree $p+1$, in the sense that
\begin{equation}
\mathcal{N}(\lambda u)=\lambda^{p+1}\mathcal{N}(u),\quad \lambda, u\in\mathbb{R}.\label{ijcm22}
\end{equation}
Thus, differentiating (\ref{ijcm22}) with respect to $\lambda$ and evaluating at $\lambda=1$, we have
\begin{equation}
\mathcal{N}^{\prime}(u)u=(p+1)\mathcal{N}(u),\quad u\in\mathbb{R}.\label{ijcm23}
\end{equation}
On the other hand, the operator $\mathcal{L}$ is invertible for $c>0$ and if $\phi=\phi_{c}$ satisfies (\ref{ijcm21}) then, using (\ref{ijcm23}), we have
\begin{equation*}
\mathcal{L}^{-1}\mathcal{N}^{\prime}(\phi_{c})\phi_{c}=(p+1)
\mathcal{L}^{-1}\mathcal{N}(\phi_{c})=(p+1)\phi_{c}.
\end{equation*}
This means that $\phi_{c}$ is an eigenfunction of the iteration operator $\mathcal{L}^{-1}\mathcal{N}^{\prime}(\phi_{c})$ with eigenvalue $p+1>1$. Therefore, for a given initial $\phi_{0}$, the classical fixed point iteration
\begin{equation*}
\mathcal{L}\phi_{n+1}=\mathcal{N}(\phi_{n}),\quad n=0,1,\ldots
\end{equation*}
will not be, in general, convergent. An alternative iterative method of fixed-point type is the Petviashvili method, \cite{Pet1976}, which is formulated as
\begin{subequations}\label{ijcm24}
\begin{equation}
m(\phi_{n})=\frac{(\mathcal{L}\phi_{n},\phi_{n})}
{(\mathcal{N}(\phi_{n}),\phi_{n})},\label{ijcm24a}
\end{equation}
\begin{equation}
\mathcal{L}\phi_{n+1}=m(\phi_{n})^{\epsilon}\mathcal{N}(\phi_{n}),\quad n=0,1,\ldots,\label{ijcm24b}
\end{equation}
\end{subequations}
for some parameter $\epsilon$. The term (\ref{ijcm24a}) is called the stabilizing factor. The convergence of (\ref{ijcm24}) for equations of the form (\ref{ijcm1}) was studied in \cite{PelinovskyS2004}. Pelinovsky and Stepanyants prove the convergence for $1<\epsilon<(p+2)/p$, under some hypotheses on the spectrum of the linearized operator $\mathcal{L}-\mathcal{N}^{\prime}(\phi_{c})$ at the profile $\phi_{c}$, with the fastest rate of convergence given by $\epsilon^{*}=(p+1)/p$. The inclusion of the stabilizing factor modifies the spectrum of $\mathcal{L}^{-1}\mathcal{N}^{\prime}(\phi_{c})$ in such a way that the eigenvalue $\lambda=p+1$ becomes, for the values of $\epsilon$ considered, an eigenvalue of the iteration operator of (\ref{ijcm24}) at $\phi_{c}$ with magnitude below one (and which is equals zero in the case of choosing $\epsilon=\epsilon^{*}$). The rest of the spectrum is preserved, see \cite{AlvarezD2014} and references therein.

The implementation of (\ref{ijcm24}) is typically carried out ´by Fourier pseudospectral approximation, \cite{AlvarezD2014,AlvarezD2016}. Let us consider the periodic problem of (\ref{ijcm21}) on a sufficiently long interval $(-l,l)$. This is discretized by a uniform grid $x_{j}=-l+jh, j=0,\ldots,N-1, h=2l/N$, in such a way that (\ref{ijcm21}) can be approximated by the discrete system
\begin{subequations}\label{ijcm25}
\begin{equation}
\mathcal{L}_{h}=
\mathcal{N}_{h}(\phi_{h}),\label{ijcm25a}
\end{equation}
\begin{equation}
\mathcal{L}_{h}=\gamma c I_{N}+D_{N}^{\alpha},\quad
\mathcal{N}_{h}(\phi_{h})=\frac{\phi_{h}.^{p+1}}{p+1},\label{ijcm25b}
\end{equation}
\end{subequations}
where $\phi_{h}$ is a $N$-vector approximation to the profile $\phi$ at the collocation points $x_{j}$, $I_{N}$ is the $N\times N$ identity matrix and $D_{N}^{\alpha}$ defined as the $N\times N$ matrix
\begin{equation*}
D_{N}^{\alpha}=F_{N}^{-1}\Lambda_{N}^{\alpha}F_{N},
\end{equation*}
with $F_{N}$ the discrete Fourier transform matrix on $\mathbb{C}^{N}$ and $\Lambda_{N}^{\alpha}$ the $N\times N$ diagonal matrix with diagonal entries of the form $|k\pi/l|^{\alpha}, k=0,\ldots,N-1$, \cite{Canutohqz,Boyd2000}. The dot in (\ref{ijcm25b}) stands for the Hadamard product. Finally, system (\ref{ijcm25a}) is implemented for the discrete Fourier coefficients, $\widehat{\phi}_{h}=(1/N)F_{N}\phi_{h}$, of $\phi_{h}$ and the resulting algebraic system is iteratively solved with the discrete version of (\ref{ijcm24})
\begin{subequations}\label{ijcm26}
\begin{equation}
m(\phi_{n})=\frac{\sum_{k=0}^{N-1}\left(c+|k\pi/l|^{\alpha}\right)|\widehat{\phi_{n}}(k)|^{2}}
{\sum_{k=0}^{N-1}\widehat{\mathcal{N}(\phi_{n})}(k)\overline{\widehat{\phi_{n}}(k)}},\label{ijcm26a}
\end{equation}
\begin{equation}
\widehat{\phi_{n+1}}(k)=m(\phi_{n})^{\epsilon}\frac{\widehat{\mathcal{N}(\phi_{n})}(k)}
{c+|k\pi/l|^{\alpha}},\label{ijcm26b}
\end{equation}
\end{subequations}
for $k=0,\ldots,N-1, n=0,1,\ldots$ and where $\widehat{\phi_{n}}=(\widehat{\phi_{n}}(0),\ldots,\widehat{\phi_{n}}(N-1))^{T}$.

A final point concerns the way how the iteration is controlled. This is done by using three strategies:
\begin{itemize}
\item Since in the case of convergence, the sequence of the stabilizing factors $m_{n}:=m(\phi_{n})$, must go to one, see (\ref{ijcm26a}), then a first control is given by the differences 
\begin{equation}
|1-m_{n}|, \quad n=0,1,\ldots\label{ijcm27a}
\end{equation}
\item A second group of control parameters is given by the sequence of the Euclidean errors between two consecutive iterations
\begin{equation}
ERROR_{c}(n)=||\phi_{n}-\phi_{n-1}||,\quad n=0,1,\ldots\label{ijcm27b}
\end{equation}
\item Finally, the sequence of the residual errors (also in Euclidean norm) 
\begin{equation}
RES(n)=||\mathcal{L}_{h}\phi_{n}-\mathcal{N}_{h}(\phi_{n})||,\quad n=0,1,\ldots,
\label{ijcm27c}
\end{equation}
is also considered.
\end{itemize}
Thus, the iteration is run up to one of (\ref{ijcm27a})-(\ref{ijcm27c}) is below a fixed tolerance parameter $tol$ which, for the experiments below, has been taken as $10^{-10}$.

\subsection{Acceleration techniques}
As mentioned in the Introduction, the loss of smoothness at the boundary due to the periodic approximation to the equations for the profiles implies an algebraic decrease of the modulus of the Fourier coefficients and, consequently, a slow iteration. In order to improve the velocity of convergence, our proposal here is including some acceleration method in the iterative process, \cite{brezinski2,sidifs}. In this sense, the so-called Vector Extrapolation Methods (VEM), \cite{sidi,smithfs,jbilous}, introduce a final stage of extrapolation at the end of each iteration of (\ref{ijcm26}). This is usually carried out in a cycling way: from the last iterate $\psi_{0}=\phi_{n}$ at stage $n$, a number $mw$ (called width of extrapolation) of iterations $\psi_{1},\ldots,\psi_{mw}$ of (\ref{ijcm26}) is computed, and the next iteration $\phi_{n+1}$ is derived as a suitable extrapolation formula from $\psi_{0},\ldots,\psi_{mw}$, see \cite{sidifs} for details. The coefficients of the extrapolation are functions of previous steps of the iteration and their derivation makes use of different criteria, leading to different methods. The one considered for the experiments  in this paper is the minimal polynomial extrapolation (MPE), a polynomial method which computes the coefficients by setting orthogonality conditions on the generalized residual, see \cite{sidi,jbilous}.


The procedure (\ref{ijcm26}), accelerated with MPE is now illustrated by the following numerical results.

\section{Numerical experiments}
\label{sec3}
\subsection{Some experiments of validation}
\begin{figure}[!htbp]
\centering
\includegraphics[width=0.7\textwidth]{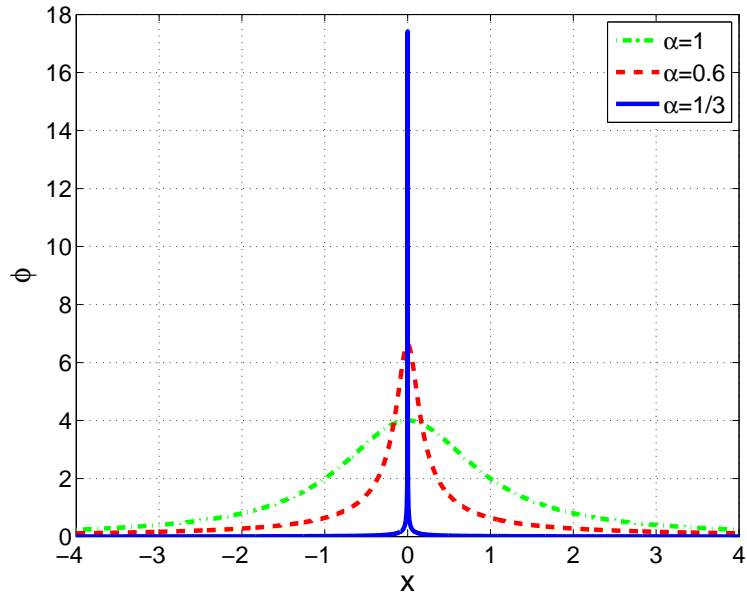}
\caption{Computed solitary wave profiles of (\ref{ijcm1}) with $p=1, c=1$ and several values of $\alpha$.}
\label{ijcm_f1}
\end{figure}

\begin{figure}[!htbp]
\centering
\subfigure[]
{\includegraphics[width=5.8cm,height=5.8cm]{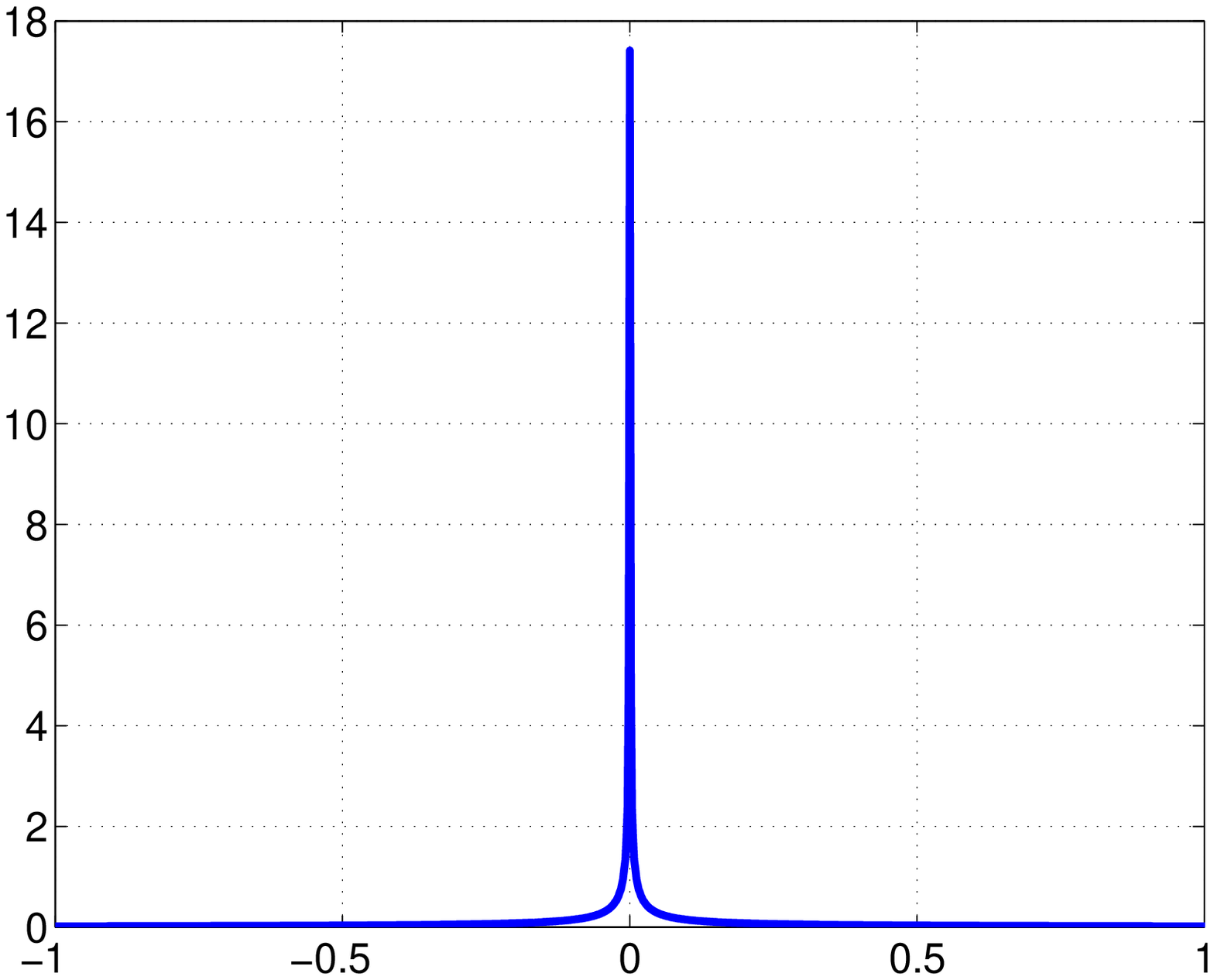}}
\subfigure[]
{\includegraphics[width=5.8cm,height=5.8cm]{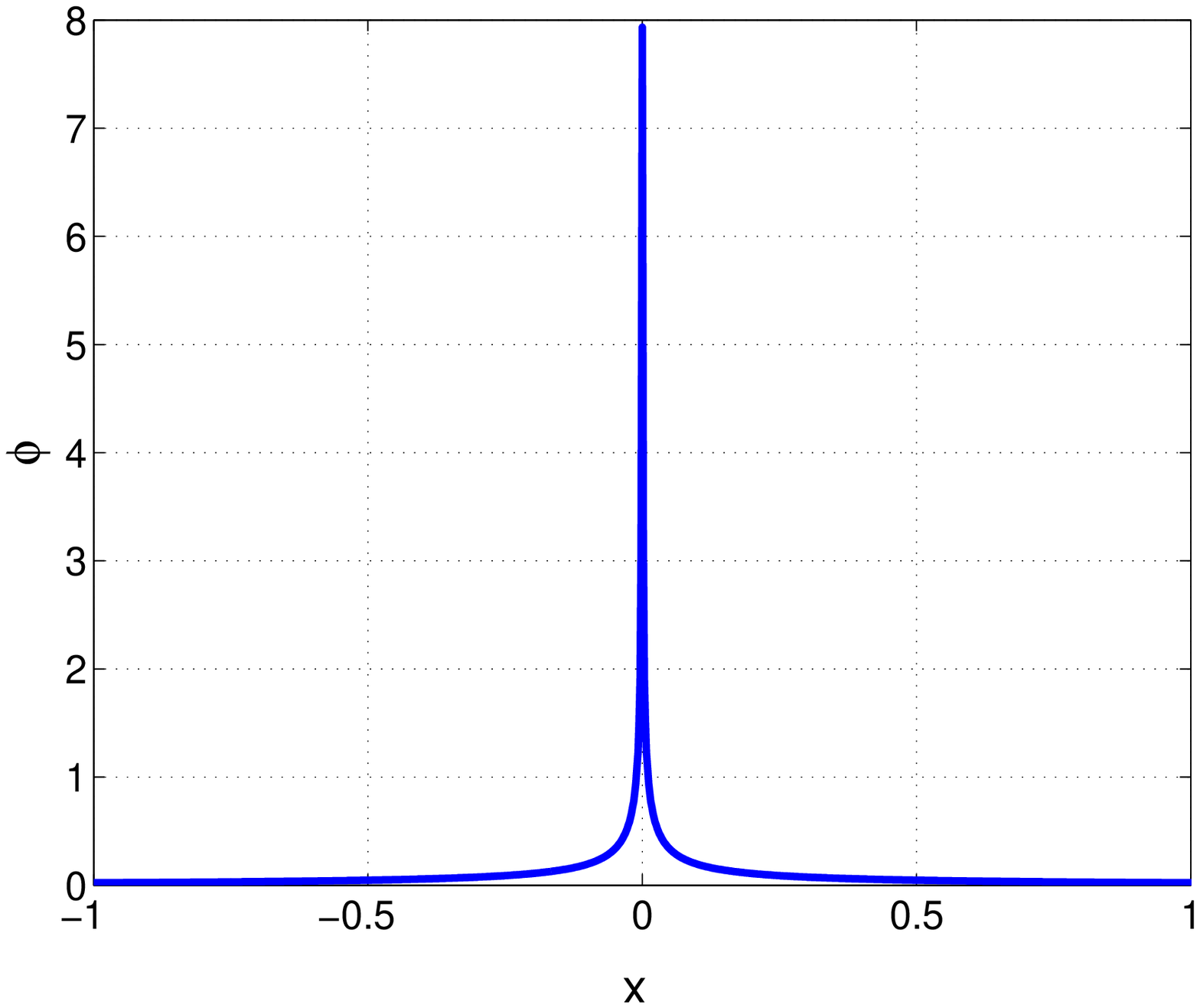}}
\subfigure[]
{\includegraphics[width=5.8cm,height=5.8cm]{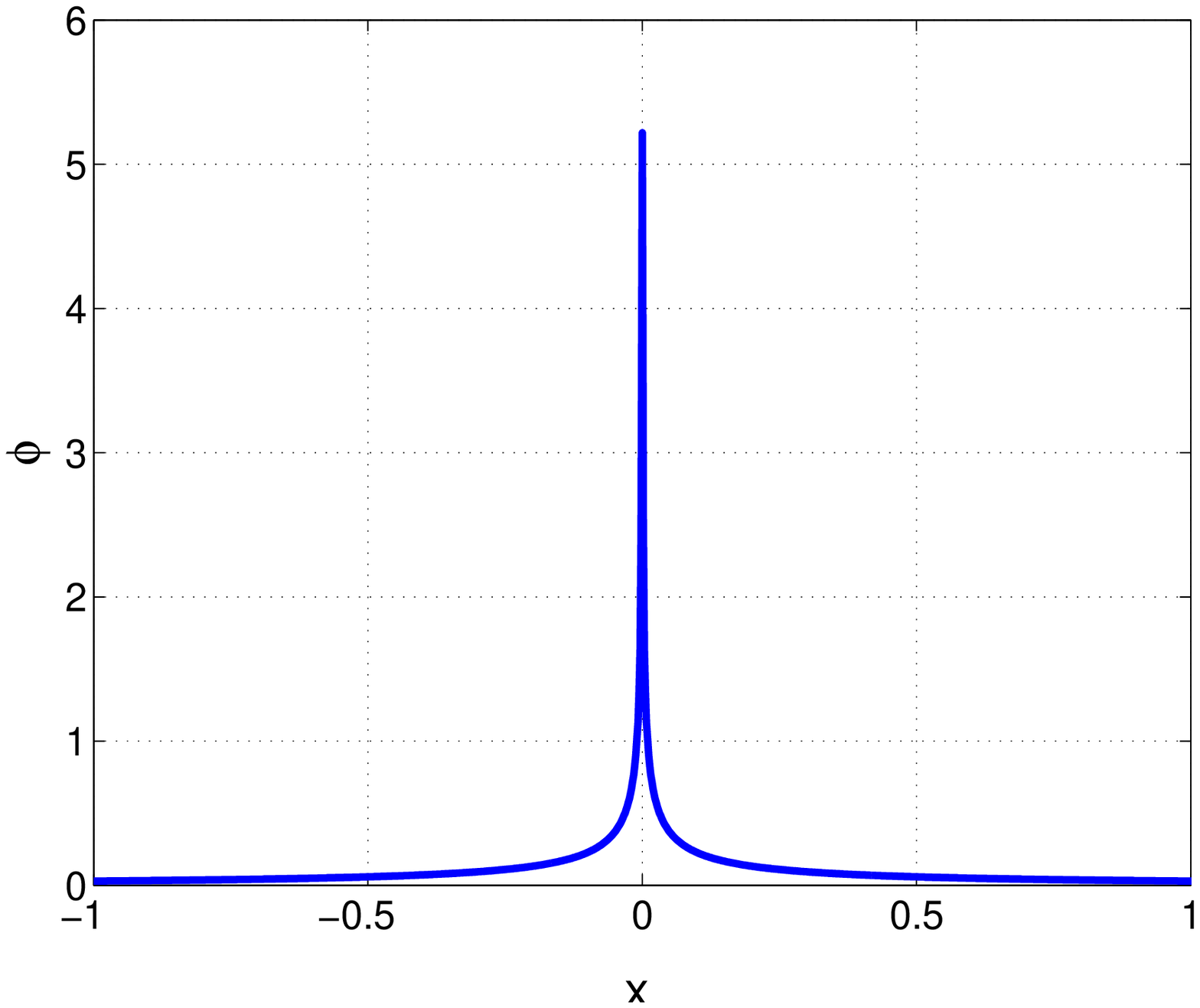}}
\subfigure[]
{\includegraphics[width=5.8cm,height=5.8cm]{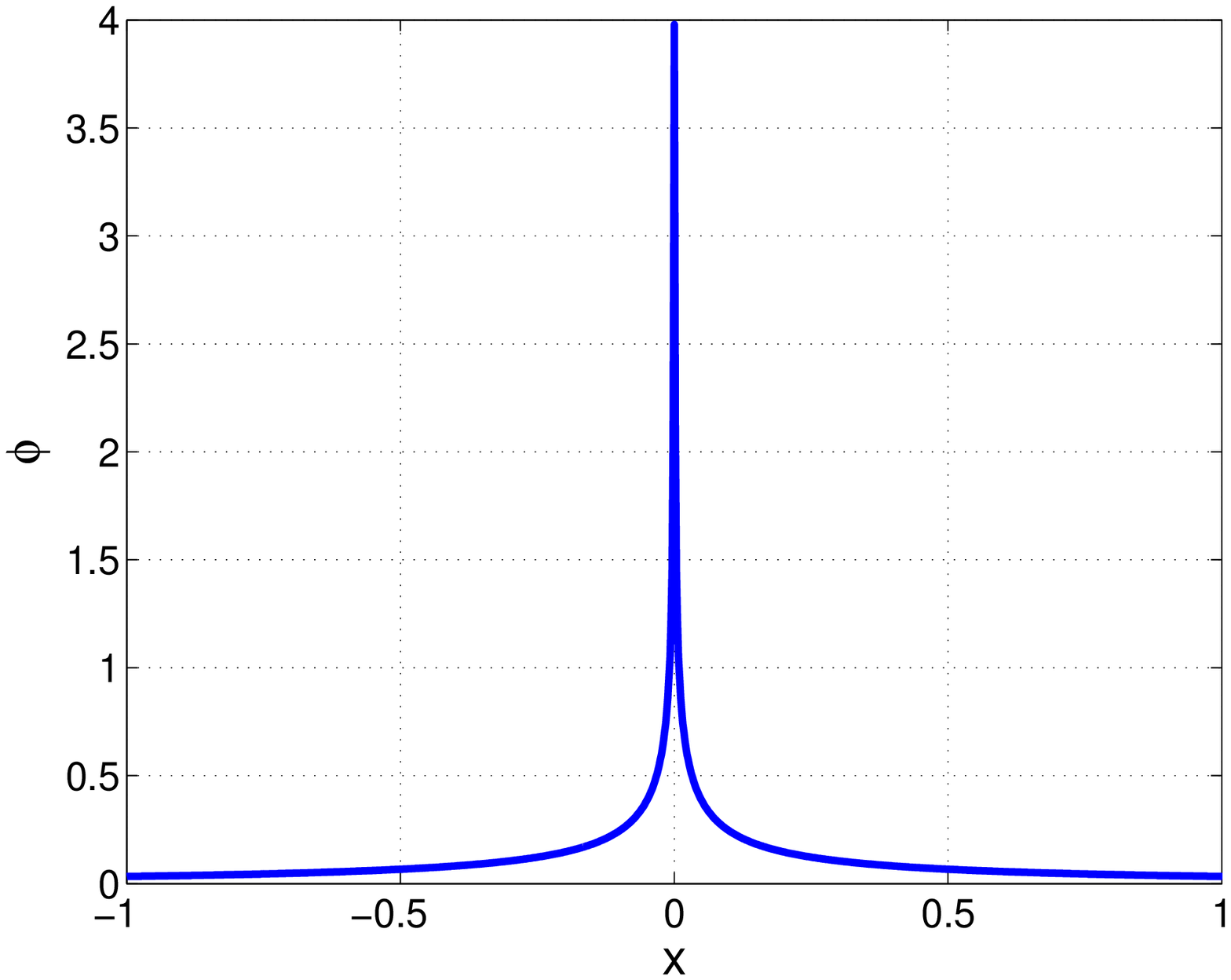}}
\caption{Computed solitary wave profiles of (\ref{ijcm1}) with $c=1$ and limiting value $\alpha_{l}=p/(p+2)$ of $\alpha$ for (a) $p=1$, (b) $p=1$, (c) $p=3$ and (d) $p=4$.}
\label{ijcm_f1b}
\end{figure}

\begin{figure}[!htbp]
\centering
\includegraphics[width=0.7\textwidth]{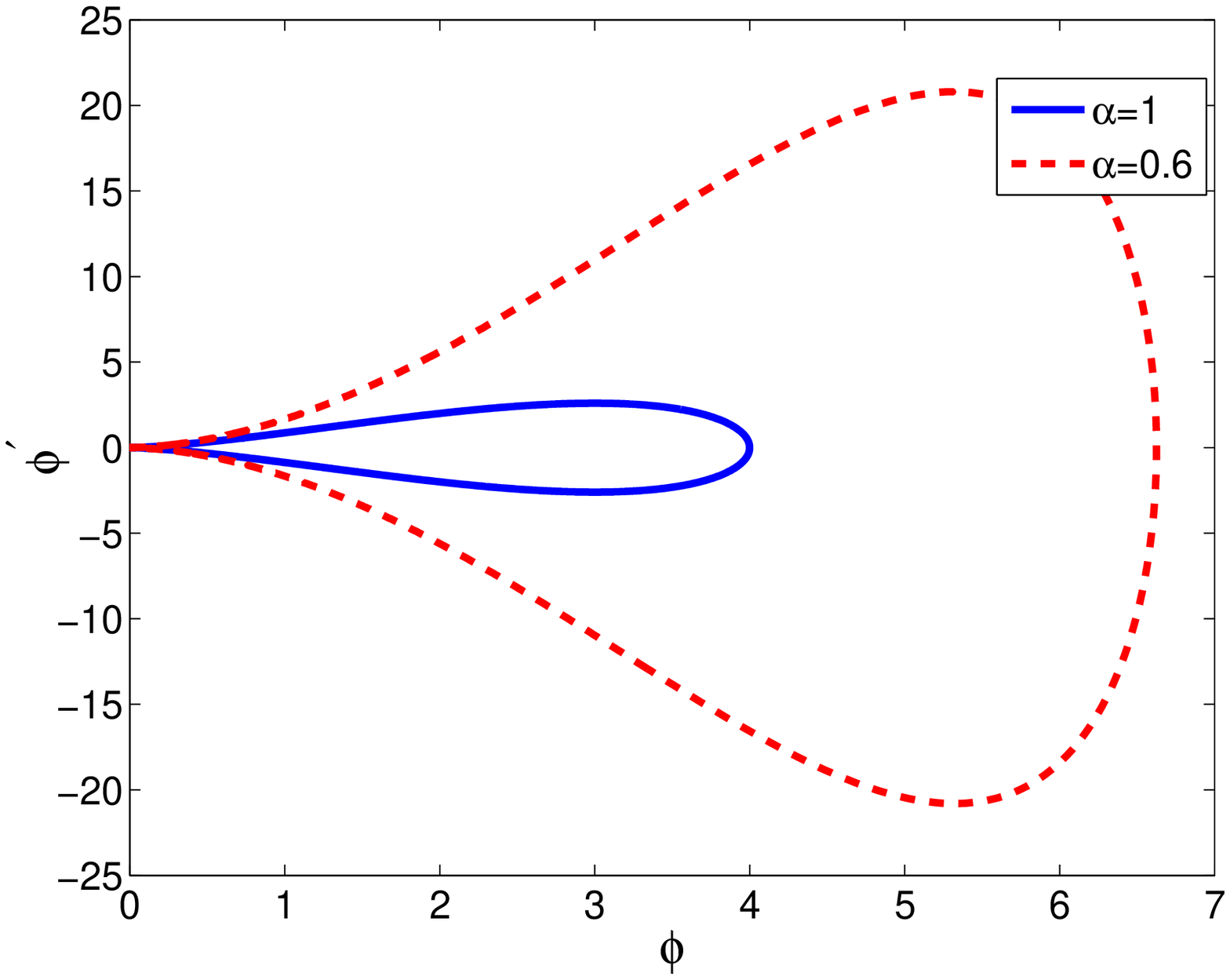}
\caption{Phase portraits of two of the numerical profiles shown in Figure \ref{ijcm_f1}.}
\label{ijcm_f2}
\end{figure}

\begin{figure}[!htbp]
\centering
\subfigure[]
{\includegraphics[width=5.8cm,height=5.8cm]{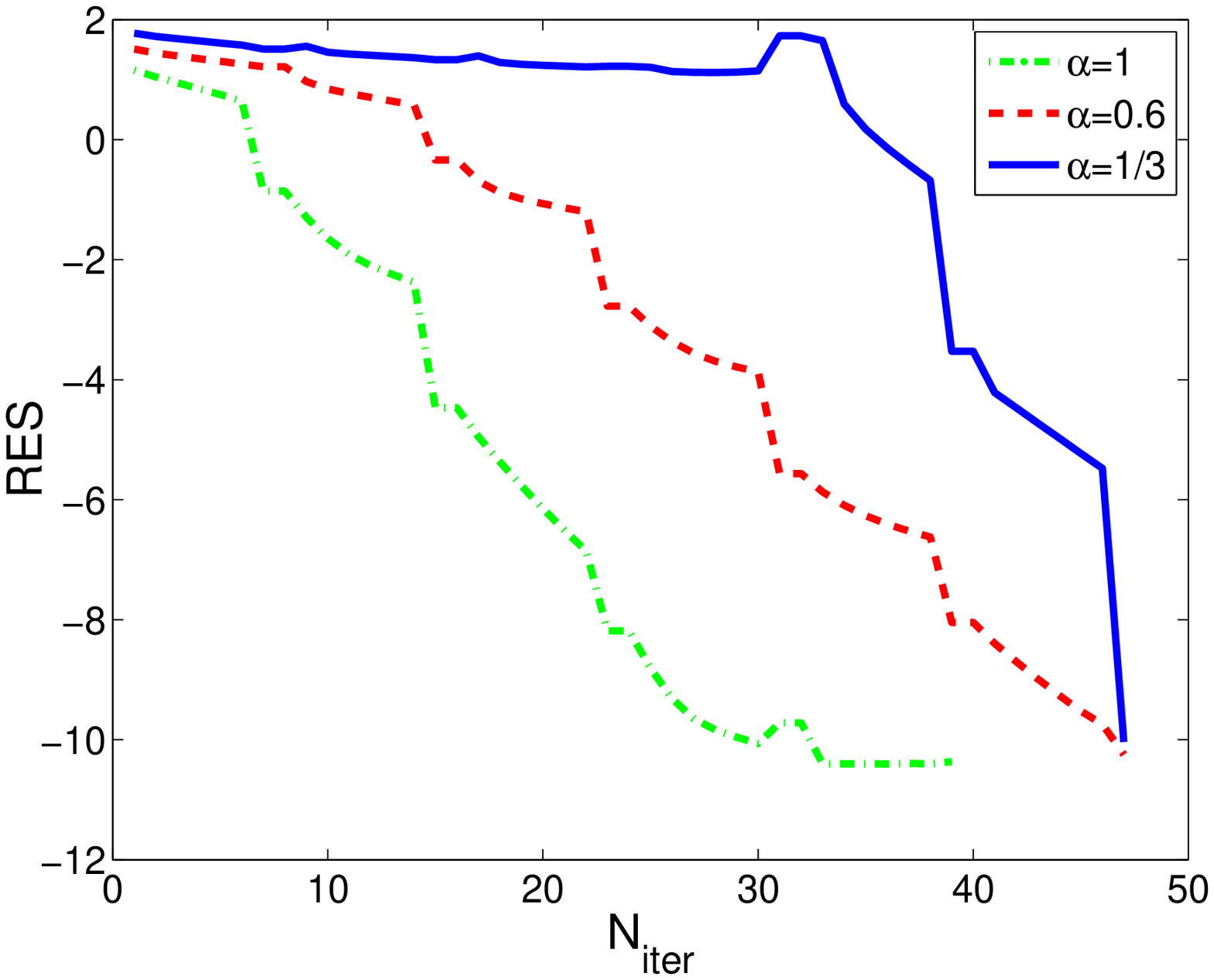}}
\subfigure[]
{\includegraphics[width=5.8cm,height=5.8cm]{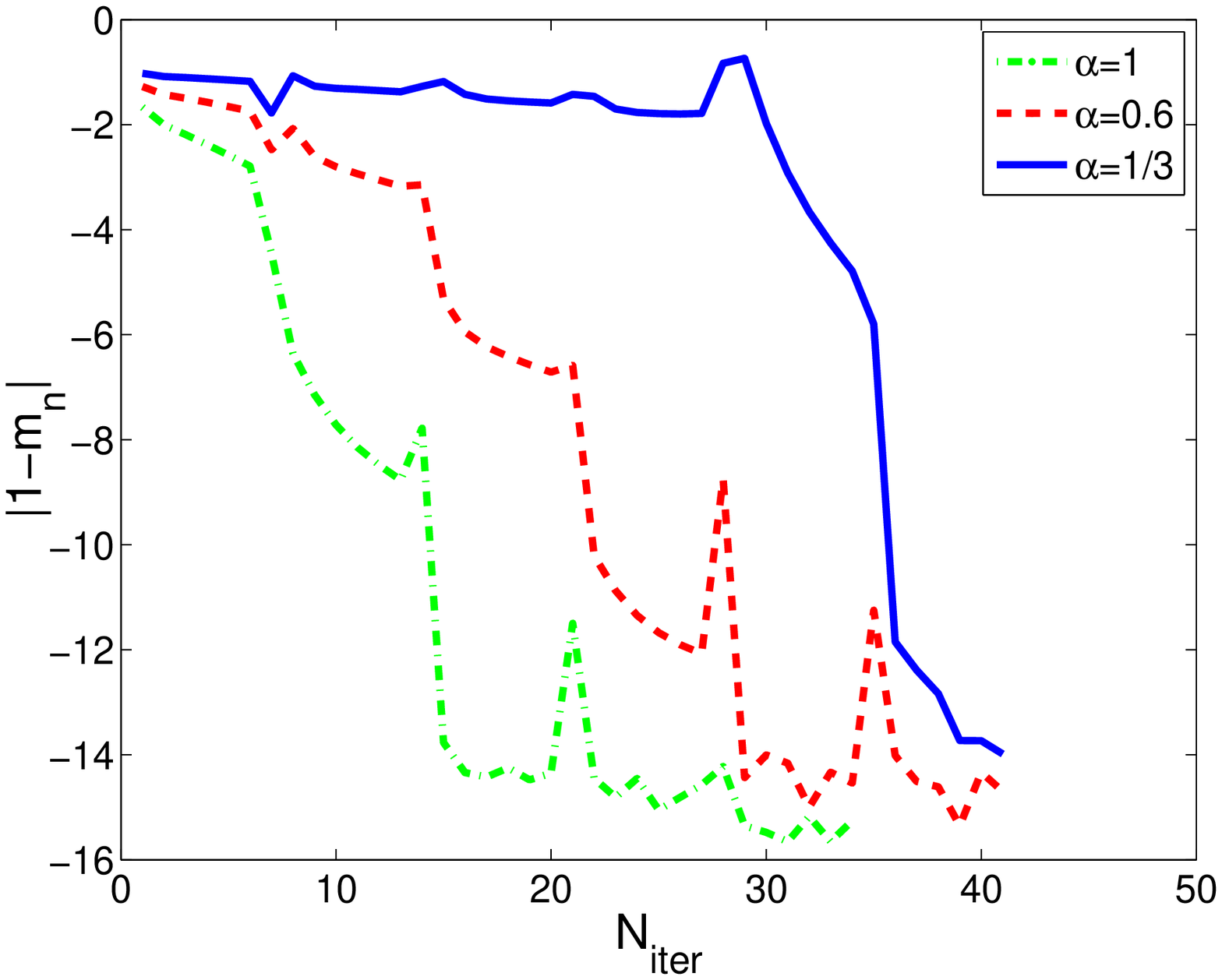}}
\caption{Generation of solitary waves of (\ref{ijcm1}) with $c=1$, $p=1$. (a) Residual errors (\ref{ijcm27c}) vs. number of iterations. (b) Stabilizing factor errors  (\ref{ijcm27a}) vs. number of iterations, both in semi-log scale.}
\label{ijcm_f3}
\end{figure}

Some experiments on the performance of the iterative scheme are first presented. Figure \ref{ijcm_f1} shows the form of approximate solitary wave profiles of (\ref{ijcm4}) with $p=1, c=1, mw=6$ and several values of $\alpha$. The initial iteration is a squared hyperbolic secant and for the implementation, an interval with $l=2048$ and $N=2^{18}$ Fourier modes were taken. As it is known (at least for $p=1$, \cite{KleinS2015}), the more peaked the profile the smaller $\alpha$ is. (This behaviour is independent of the nonlinearity parameter, although as $p$ increases, the amplitude of the corresponding profile has been observed to decrease, see figures \ref{ijcm_f6} and \ref{ijcm_f7} below.) In all the computations, the minimum value of the profile is below $10^{-5}$. 
From Figure \ref{ijcm_f1}, note that the limiting case $\alpha_{l}=1/3$, for $p=1$, of $\alpha$ looks to be computable. This also holds for any $p>1$, for which $\alpha_{l}=p/(p+2)$, see Figure \ref{ijcm_f1b}.
The algebraic decay at infinity of some of the profiles can be observed in Figure \ref{ijcm_f2}, which displays the corresponding phase portraits and where the derivative has been computed by using pseudospetral differentiation, see \cite{Canutohqz,Boyd2000}.
In order to check the accuracy of the computed profiles, several experiments are made. Figure \ref{ijcm_f3} displays the behaviour of the residual error and the stabilizing factor as function of the number of iterations and for each of the waves computed in Figure \ref{ijcm_f1}. The results confirm the convergence of the sequences (\ref{ijcm27a}), (\ref{ijcm27c}) and consequently of the iteration. (The second control sequence (\ref{ijcm27b}) behaves even better and the corresponding results are not shown.) As far as the performance is concerned, note that in all the cases, the tolerance $tol=10^{-10}$ is attained in less than $50$ iterations.
A second experiment to check the accuracy is as follows. The computed profiles were taken as initial condition of a numerical method to approximate the periodic initial value problem associated to (\ref{ijcm1}). The numerical scheme is based on a pseudospectral Fourier discretization in space and a fourth order, diagonally implicit, Runge-Kutta composition method, described in \cite{dFrutosS1992} (see also references therein) as time integrator. The scheme has relevant geometric properties, \cite{HairerLW2004}, and has been shown to be efficient in nonlinear wave problems, \cite{KalischMV2017}.

\begin{figure}[!htbp]
\centering
\subfigure[]
{\includegraphics[width=5.8cm,height=5.8cm]{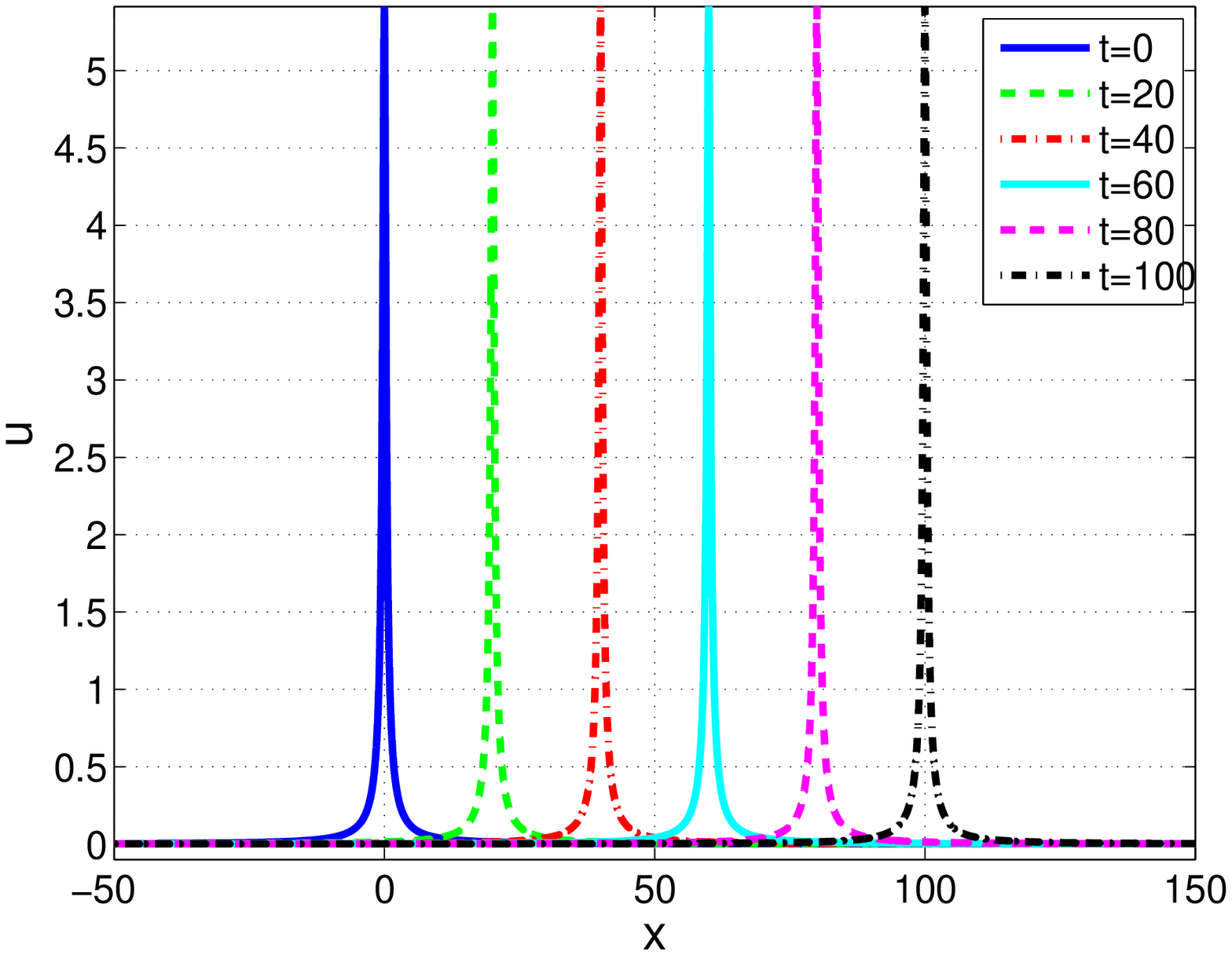}}
\subfigure[]
{\includegraphics[width=5.8cm,height=5.8cm]{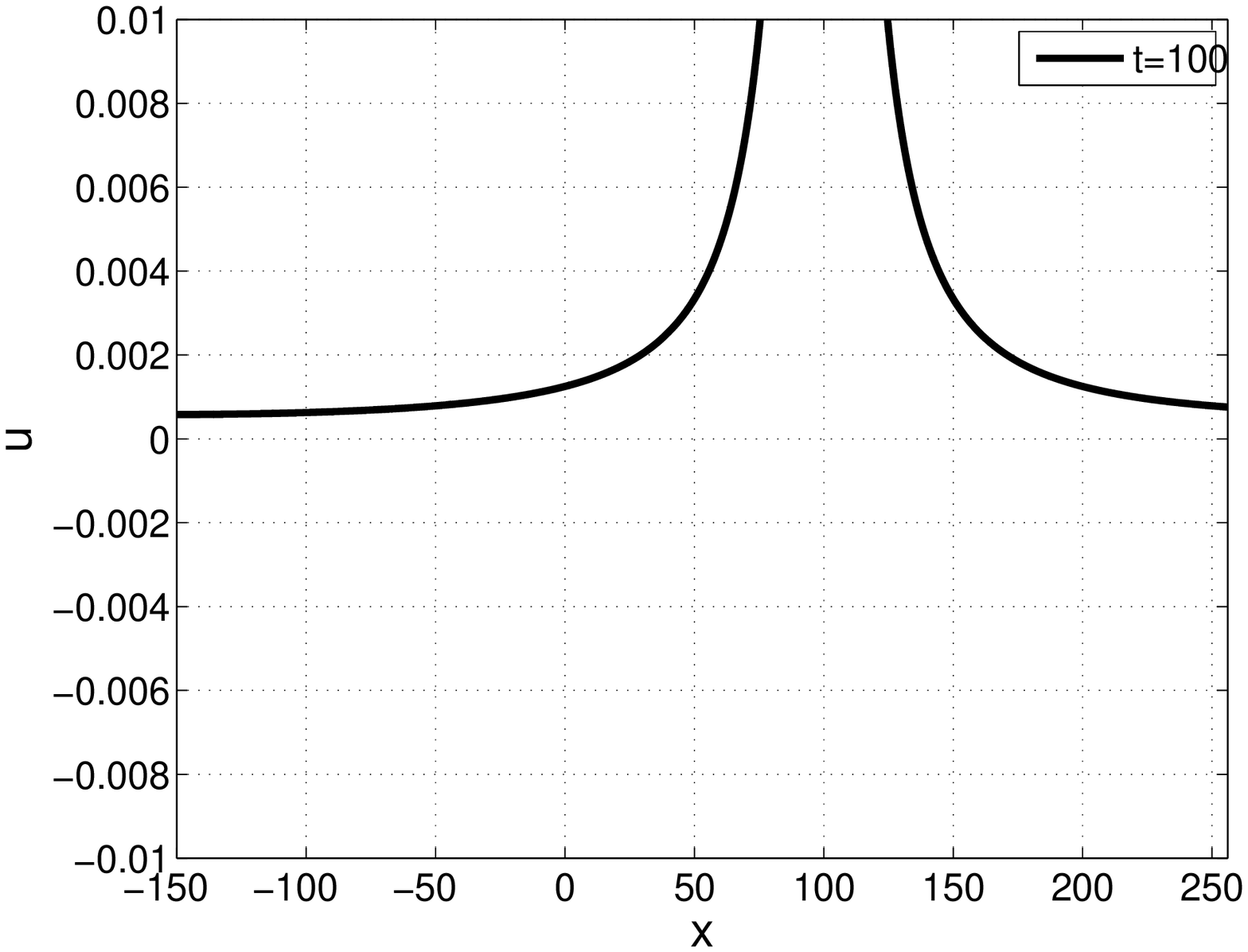}}
\subfigure[]
{\includegraphics[width=5.8cm,height=5.8cm]{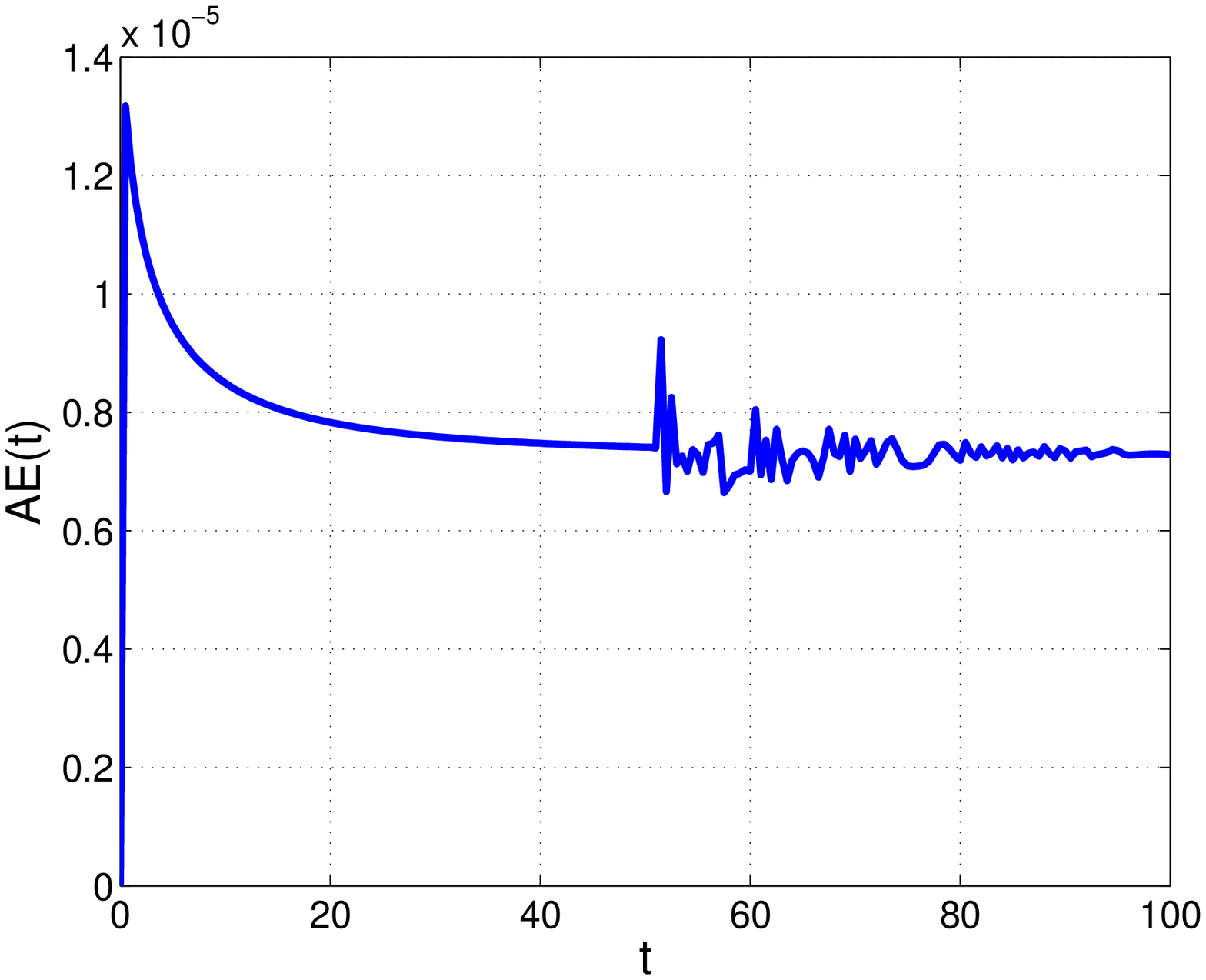}}
\subfigure[]
{\includegraphics[width=5.8cm,height=5.8cm]{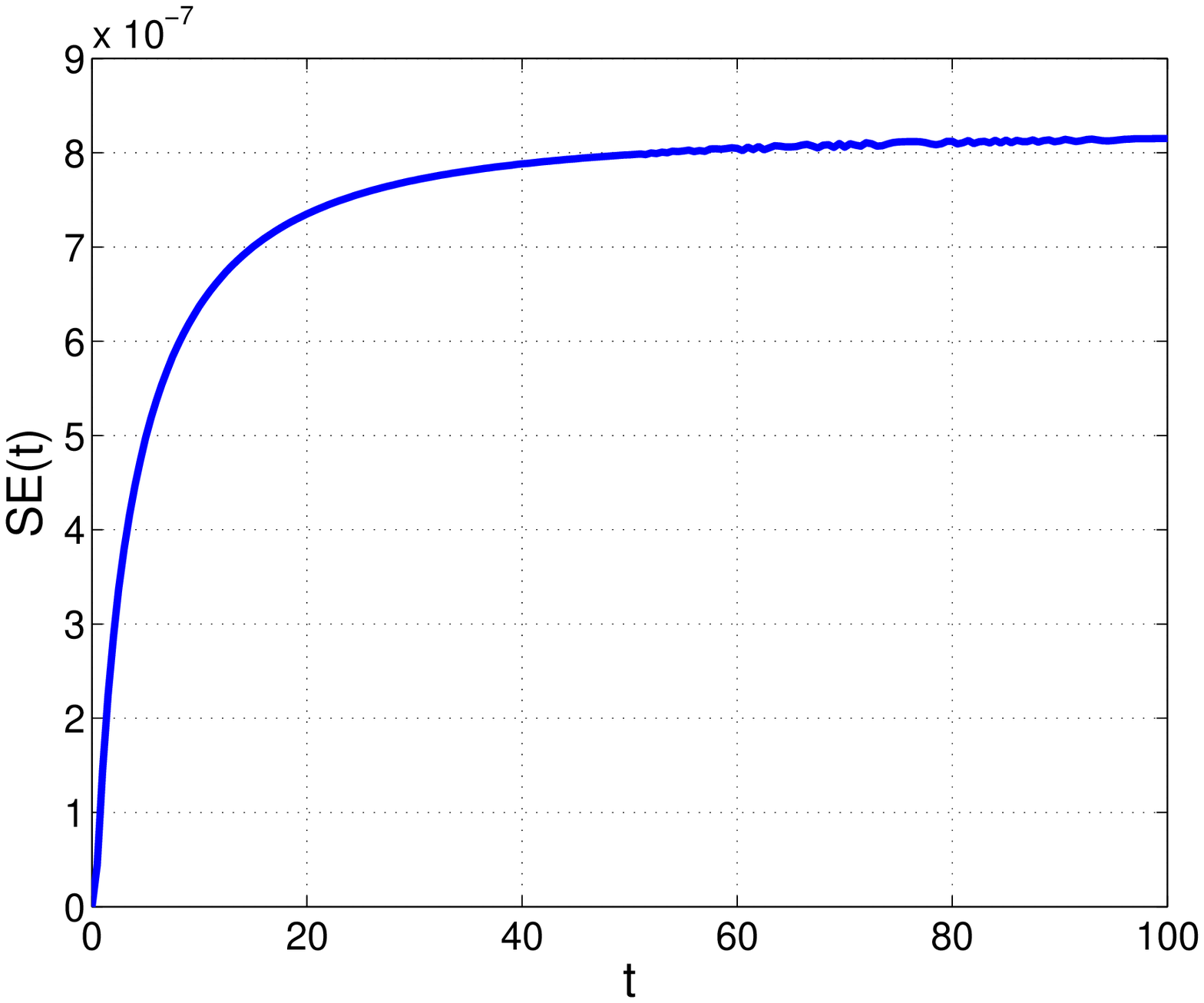}}
\caption{(a) Numerical approximation of (\ref{ijcm1}) from the profile computed with $\alpha=0.7, c=1$ and $p=1$, at several times. (b) Magnification of (a). (c) Amplitude and (d) speed  errors vs. time.}
\label{ijcm_f4}
\end{figure}

\begin{figure}[!htbp]
\centering
\includegraphics[width=0.7\textwidth]{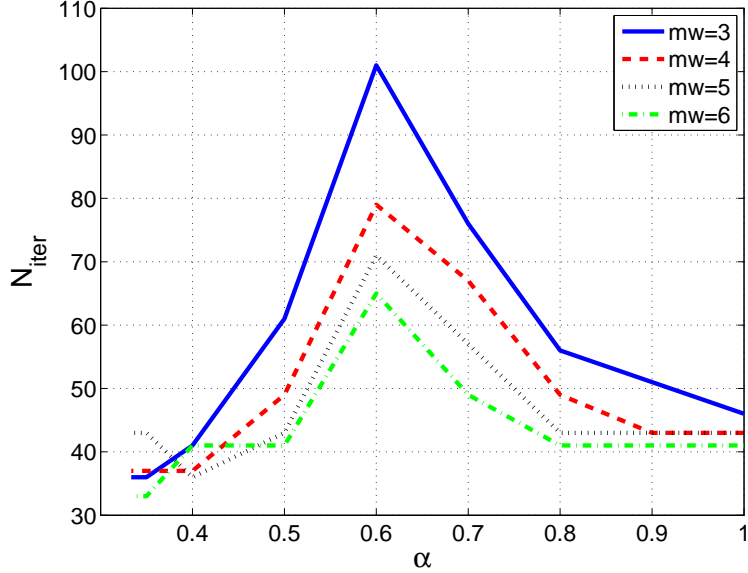}
\caption{Computed solitary wave profiles of (\ref{ijcm1}) with $p=1, c=1$. Number of iterations to attain a residual error below $tol=10^{-10}$ as function of $\alpha$ and for several values of $mw$.}
\label{ijcm_f5}
\end{figure}

The evolution of the corresponding numerical solution was monitored and, in the case of $\alpha=0.7$, is represented at several times in Figures \ref{ijcm_f4}(a),(b). They show that the initial approximate profile evolves in a solitary way without relevant disturbances, suggesting that the computed profile represents a solitary wave of (\ref{ijcm1}) with a high degree of accuracy. This is also confirmed by the evolution of the amplitude and speed of the numerical approximation during the integration, displayed in Figures \ref{ijcm_f4}(c) and (d). (The computation of amplitude and speed was made in the standard way, see e.~g. \cite{DougalisDLM2007}.)
The following experiments complement the illustration of the performance of the method. Figure \ref{ijcm_f5}
shows the number of iterations required to get a residual error below the tolerance as function of the fractional parameter $\alpha$ and for several values of the width of extrapolation $mw$, in the approximation to a solitary wave profile of (\ref{ijcm4}) with $p=1$ and $c=1$. Note that, for moderate values of $mw$, as $mw$ increases, the number of iterations decreases and in this sense the performance improves. However, observe that the parameter $mw$ cannot be fixed \emph{a priori} in this sort of nonlinear problems, \cite{sidifs,smithfs}. Also, from some value of $mw$, the number of iterations and the computational time will not improve anymore. For one of the values of $mw$ considered, the number of iterations does not increase as $\alpha$ decreases; it is maximum when $\alpha$ is close to $0.6$.
\subsection{Application to study additional properties of the solitary waves}
\begin{figure}[!htbp]
\centering
\subfigure[]
{\includegraphics[width=5.8cm,height=5.8cm]{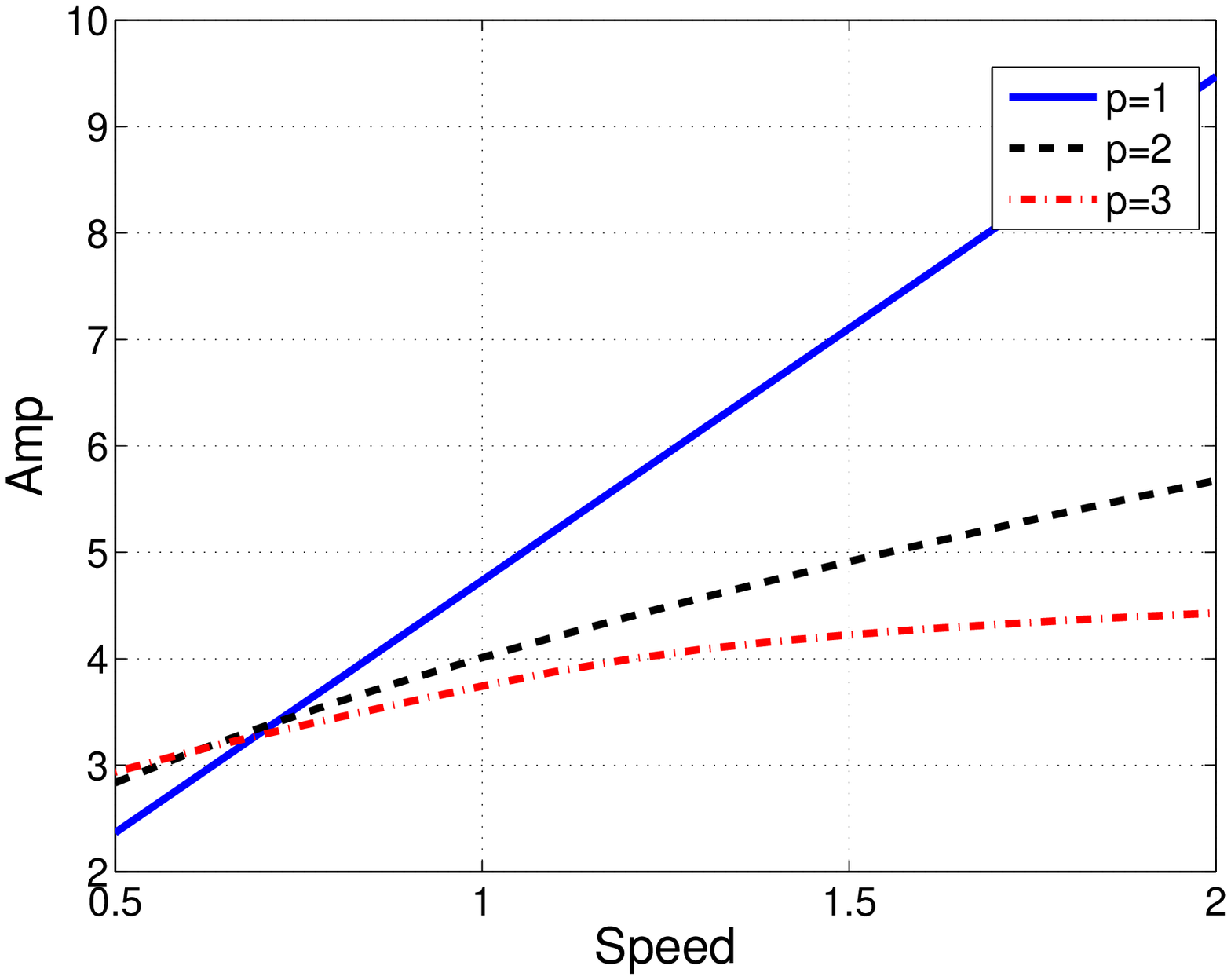}}
\subfigure[]
{\includegraphics[width=5.8cm,height=5.8cm]{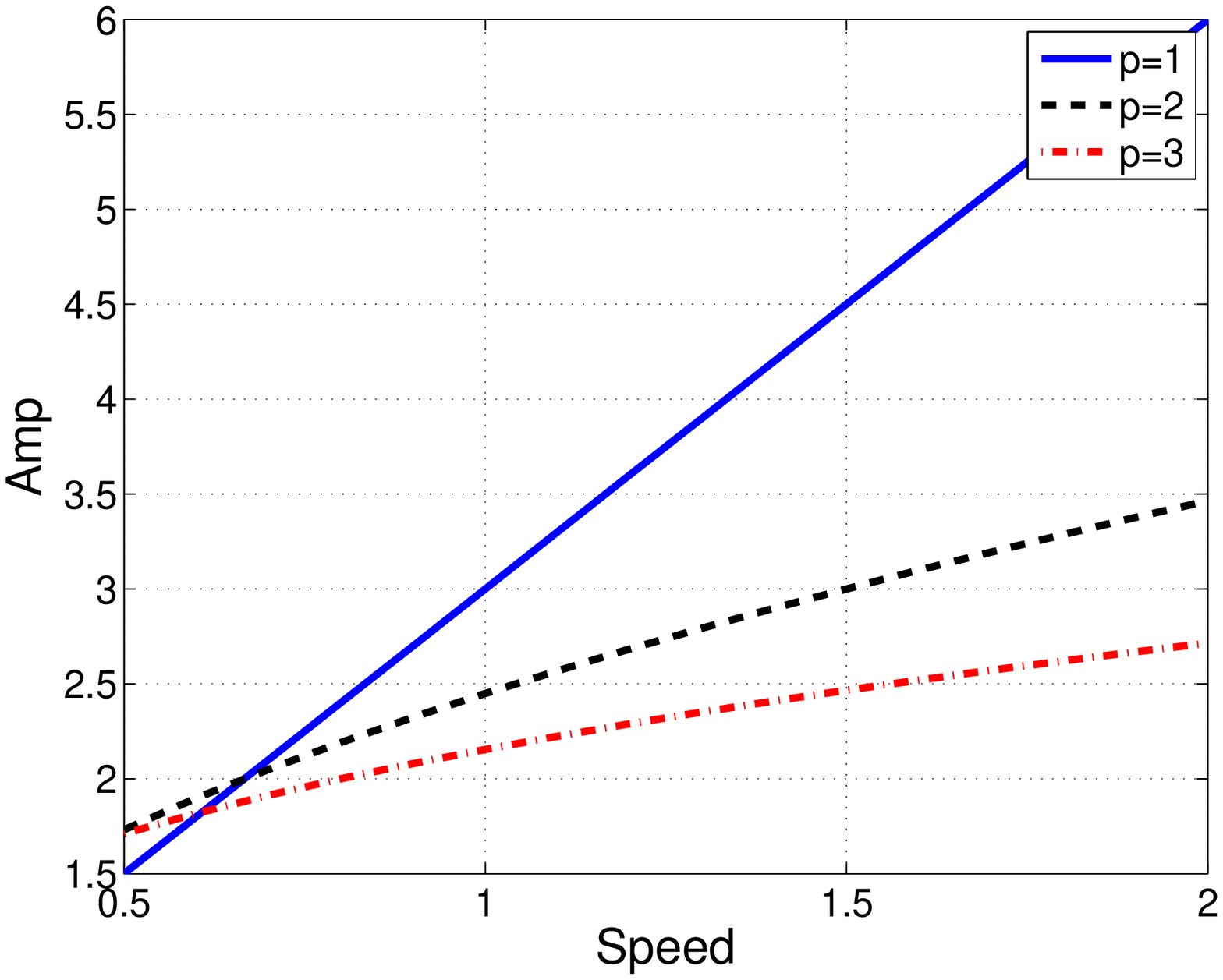}}
\caption{Speed-amplitude relation of computed solitary waves of (\ref{ijcm1}) for several values of $p$. (a) $\alpha=0.8$. (b) $\alpha=2$ (KdV case).}
\label{ijcm_f6}
\end{figure}

The final group of experiments is concerned with the model (\ref{ijcm1}). Once the accuracy and performance of the iterative method have been checked, the procedure can be used to obtain some additional information about the solitary waves. This idea is focused here on the search for the speed-amplitude relation and its dependence on the parameters $\alpha$ and $p$. Figure \ref{ijcm_f6}(a) displays this relation for $\alpha=0.8$ fixed and several values of $p$. The amplitude is an increasing function of the speed but the increment depends on $p$. When $\alpha=2$ (gKdV case) the solitary waves are known explicitly and the amplitude has the exact formula, \cite{BonaDKM1995}
\begin{equation}
Amp=\left(\frac{c}{2}(p+1)(p+2)\right)^{\frac{1}{p}},\label{ijcm_31}
\end{equation}
see Figure \ref{ijcm_f6}(b). 

\begin{figure}[!htbp]
\centering
\subfigure[]
{\includegraphics[width=5.8cm,height=5.8cm]{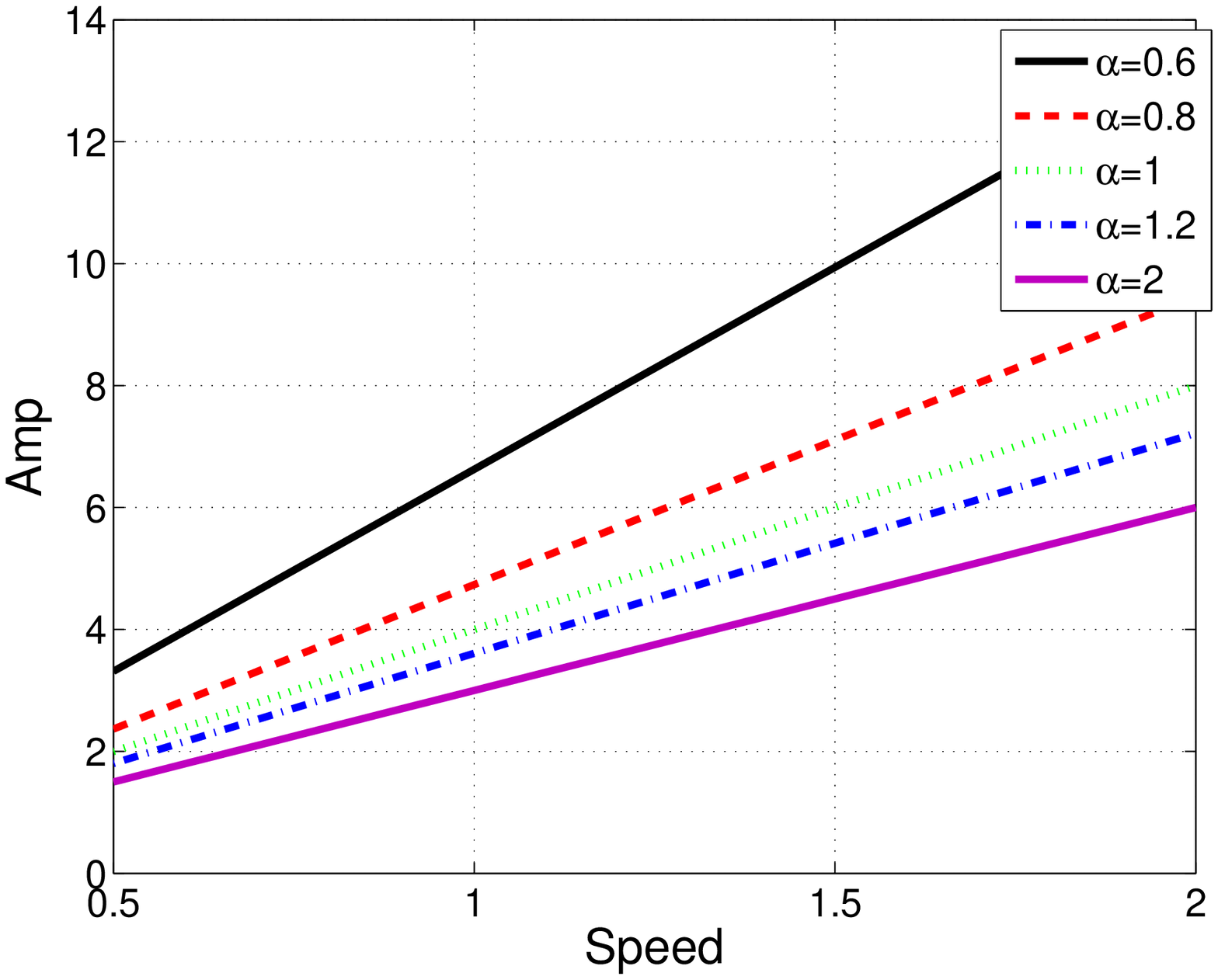}}
\subfigure[]
{\includegraphics[width=5.8cm,height=5.8cm]{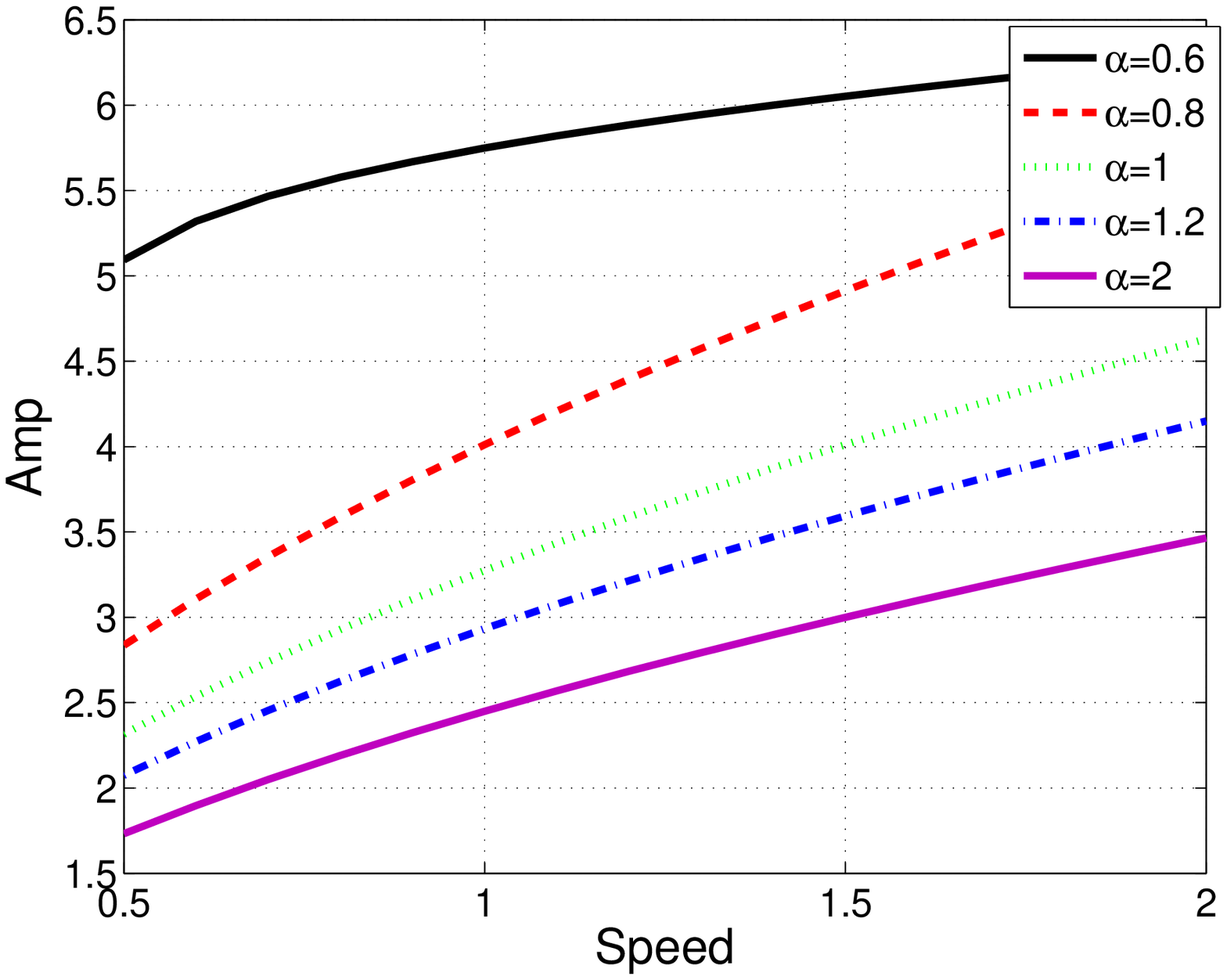}}
\subfigure[]
{\includegraphics[width=5.8cm,height=5.8cm]{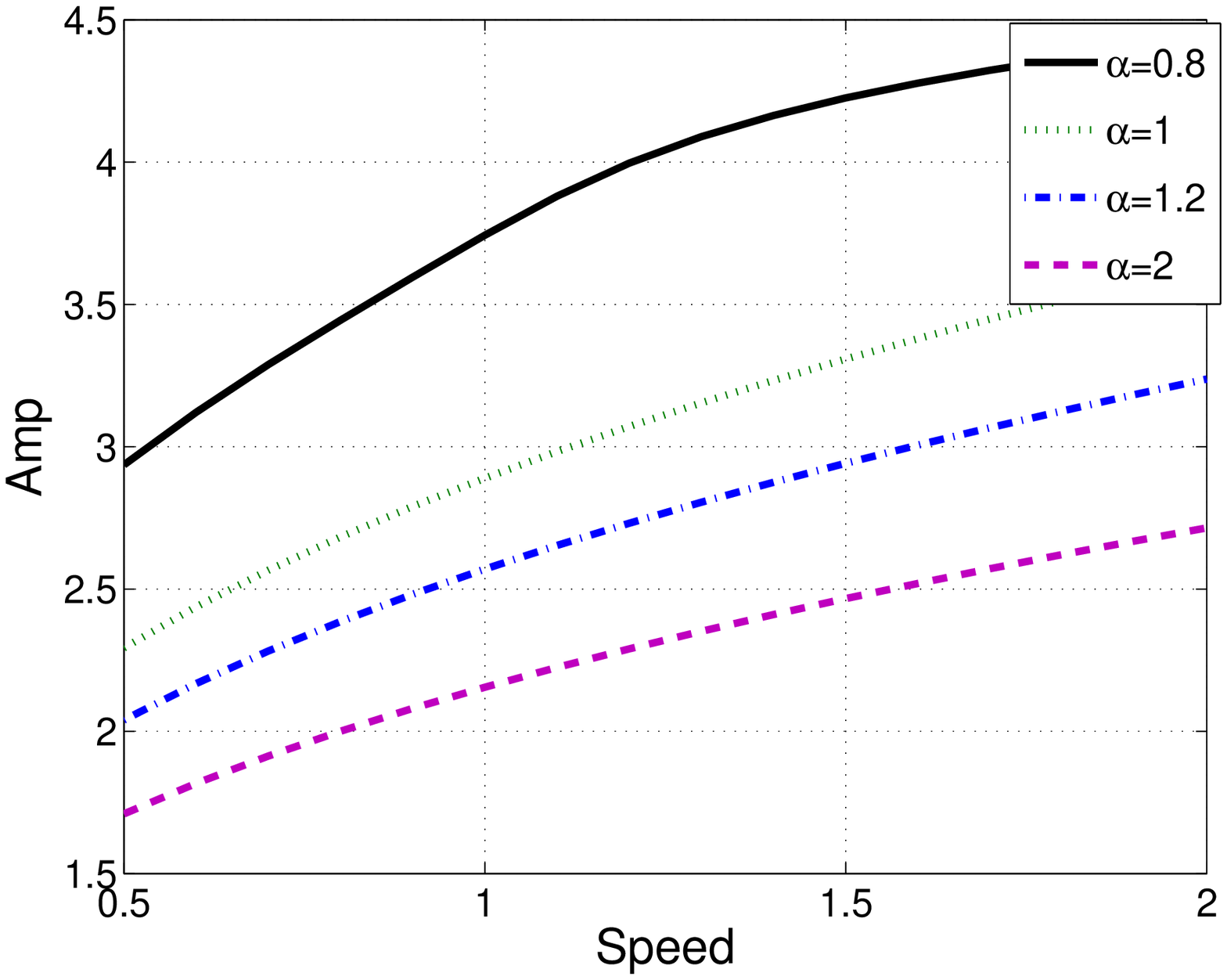}}
\caption{Speed-amplitude relation of computed solitary waves of (\ref{ijcm1}) for several values of $\alpha$. (a) $p=1$. (b) $p=2$. (c) $p=3$.}
\label{ijcm_f7}
\end{figure}

A comparison of both figures suggests to consider the experiment of studying the speed-amplitude relation for a fixed value of $p$ and as function of $\alpha$. This is illustrated by Figure \ref{ijcm_f7}. Note that for a fixed value of $c$ and $\alpha$, increasing $p$ typically leads to a solitary wave of smaller amplitude. When $c$ and $p$ are fixed, the taller the wave the smaller value of $\alpha$ is (cf. Figures \ref{ijcm_f1} and \ref{ijcm_f6}). 
Due to the peaked form of the profiles, the maximum is, for small values of $\alpha$, not easy to compute with a minimum of accuracy. However, the most reliable numerical results suggest that a similar relation to (\ref{ijcm_31}) holds for any $\alpha$. This is observed in Table \ref{ijcm_tav1} where, for several values of $p$ and $\alpha$, the numerical data speed-amplitude were fitted to a power function $f(x)=ax^{b}$. The accuracy of the results was guaranteed by a goodness of fit where some statistical parameters are around some fixed tolerance. Explicitly, the following statistics were used: the sum of squares due to error SSE (with a tolerance threshold of $10^{-6}$), the R-squared (which in all the cases is equals $1$) and the root mean squared error RSME (around $10^{-4}$).

\begin{table}
\tbl{Parameters of the fit $f(x)=ax^{b}$ of the computed speed-amplitude relation for different values of $\alpha$ and $p$.}
{\begin{tabular}{c|ccc} 
 & $p=1$ & $p=2$ & $p=3$ \\ \midrule
 $\alpha=0.8$ & $a=4.735$ & $a=4.012$& $a=3.702$ \\
  & $b=1$ & $b=0.5001$& $b=0.2966$ \\\midrule
 $\alpha=1.2$ & $a=3.61$ & $a=2.934$& $a=2.57$ \\
  & $b=1$ & $b=0.5$& $b=0.3333$ \\ \bottomrule
\end{tabular}}
\label{ijcm_tav1}
\end{table}
Note that while the coefficients $b$ suggest an amplitude as a power $1/p$ of the speed as in (\ref{ijcm_31}), the dispersion parameter $\alpha$ looks to affect only the coefficient $a$ of the fit.

\section*{Acknowledgement(s)}
This work was supported by  Spanish Ministerio de Econom\'{\i}a y Competitividad under the Research Grant MTM2014-54710-P.

\end{document}